\documentclass[12pt, reqno]{amsart}
\usepackage{indentfirst, amssymb, amsmath, amsthm, mathrsfs, setspace, indentfirst, enumerate,  mathrsfs, amsmath, amsthm}

\usepackage[bookmarksnumbered, colorlinks, plainpages]{hyperref}
\usepackage{mathrsfs}
\usepackage{cite}
\usepackage{graphicx}
\usepackage{float}
\newtheorem*{theoA}{Theorem A}
\newtheorem*{theoB}{Theorem B}

\newtheorem{ques}{Question}[section]

\newtheorem{theo}{Theorem}[section]
\newtheorem{lem}{Lemma}[section]

\newtheorem{defi}{Definition}[section]
\newtheorem{rem}{Remark}[section]

\newcommand{\ol}{\overline}

\newcommand{\be}{\begin{equation}}
\newcommand{\ee}{\end{equation}}
\newcommand{\beas}{\begin{eqnarray*}}
\newcommand{\eeas}{\end{eqnarray*}}
\newcommand{\bea}{\begin{eqnarray}}
\newcommand{\eea}{\end{eqnarray}}

\numberwithin{equation}{section}
\begin{document}

\title[D\MakeLowercase{ifferential Operators, Multiple Schwarz Functions, and the Bohr Radius of Stable Harmonic Maps}]
{Differential Operators, Multiple Schwarz Functions, and the Bohr Radius of Stable Harmonic Maps}
\author[S. M\MakeLowercase{ajumder}, N. S\MakeLowercase{arkar},  \MakeLowercase{and} M. B. A\MakeLowercase{hamed}]{S\MakeLowercase{ujoy} M\MakeLowercase{ajumder}, N\MakeLowercase{abadwip} S\MakeLowercase{arkar} \MakeLowercase{and} M\MakeLowercase{olla} B\MakeLowercase{asir} A\MakeLowercase{hamed}$^{*}$}

\address{Department of Mathematics, Raiganj University, Raiganj, West Bengal-733134, India.}
\email{sm05math@gmail.com, sjm@raiganjuniversity.ac.in}

\address{Department of Mathematics, Raiganj University, Raiganj, West Bengal-733134, India.}
\email{naba.iitbmath@gmail.com}

\address{Molla Basir Ahamed,
	Department of Mathematics,
	Jadavpur University,
	Kolkata-700032, West Bengal, India.}
\email{mbahamed.math@jadavpuruniversity.in}

\date{}

\renewcommand{\thefootnote}{}
\footnote{*\emph{Corresponding Author}: Molla Basir Ahamed.}

\renewcommand{\thefootnote}{\arabic{footnote}}
\setcounter{footnote}{0}

\begin{abstract}
In this paper, we study the Bohr phenomenon for differential operators $D$ and $\mathscr{D}$ of stable harmonic mappings involving multiple Schwarz functions in $\mathcal{B}_n$, using distance formulations. By constructing suitable combinations of multiple Schwarz functions, we establish sharp and improved Bohr-type inequalities for these mappings. The corresponding Bohr radii are also determined for certain subclasses of stable harmonic functions and their associated differential operators. Moreover, Bohr-Rogosinski-type inequalities are derived, which highlights the influence of multiple Schwarz functions on the geometric properties of stable harmonic mappings. All the radii are determined, and we prove that each one is the best possible.
\end{abstract}

\maketitle
\tableofcontents
\section{{\bf Introduction}}
The origin of the Bohr phenomenon lies in the seminal work by Harald Bohr
\cite{Bohr-1914} in $ 1914 $ for the analytic functions of the form $ \sum_{n=0}^{\infty}a_nz^n $ defined on the unit disk $ \mathbb{D}:=\{z\in\mathbb{C} : |z|<1\} $ with $ |f(z)|<1 $. Study of Bohr phenomenon for different classes of functions with various settings becomes a subject of great interests during past several years and an extensive research work has been done by many authors but Bohr phenomenon is not much investigated for the classes of harmonic mappings. In this paper, we are mainly interested to study Bohr phenomenon with suitable settings in order to establish certain harmonic analogue of some Bohr inequality valid for analytic functions.\vspace{1.2mm}

For a continuously differentiable complex-valued mapping $ f(z)=u(z)+iv(z) $, $ z=x+iy $, we use the common notions for its formal derivatives:
\begin{align*}
	f_{z}=\frac{1}{2}\left(f_x-if_y\right)\;\; \mbox{and}\;\; f_{\bar{z}}=\frac{1}{2}\left(f_x+if_y\right).
\end{align*}
We say that $ f $ is a harmonic mapping in a simply connected domain $ \Omega $ if $ f $ is twice continuously differentiable and satisfies the Laplacian equation $ \Delta f=4f_{z\bar{z}}=0 $ in $ \Omega $, where $ \Delta $ is the complex Laplacian operator defined by
$ \Delta={\partial^2}/{\partial x^2}+{\partial^2}/{\partial y^2}. $\vspace{1.2mm}

Methods of harmonic mappings have been applied to study and solve fluid flow problems (see \cite{Aleman-2012,Constantin-2017}). For instance, in 2012, Aleman and Constantin \cite{Aleman-2012} established a connection between harmonic mappings and ideal fluid flows. They developed an ingenious technique to solve the incompressible two-dimensional Euler equations in terms of univalent harmonic mappings (see \cite{Constantin-2017} for details).\vspace{1.2mm}

Let $ \mathcal{H}(\Omega) $ be the class of functions harmonic in $ \Omega $. It is well-known that functions in the class $ \mathcal{H}(\Omega) $ has the following representation $ f=h+\overline{g} $, where $ h $ and $ g $ both are analytic functions in $ \Omega $. The famous Lewy's theorem \cite{Lew-BAMS-1936} in $ 1936 $ states that a harmonic mapping $ f=h+\overline{g} $ is locally univalent on $ \Omega $ if, and only if, the determinant $ |J_f| $ of its Jacobian matrix $ J_f $ does not vanish on $ \Omega $, where
\begin{equation*}
	|J_f|:=|f_{z}|^2-|f_{\bar{z}}|^2==|h^{\prime}(z)|^2-|g^{\prime}(z)|^2\neq 0.
\end{equation*}
In view of this result, a locally univalent harmonic mapping  is sense-preserving if $ |J_f(z)|>0 $ and sense-reversing if $|J_{f}(z)|<0$ in $\Omega$. For detailed information about the harmonic mappings, we refer the reader to \cite{Clunie1984,Duren2004}. In \cite{KPSM1018}, Kayumov \emph{et al.} first established the harmonic extension of the classical Bohr theorem, since then investigating on the Bohr-type inequalities for certain class of harmonic mappings becomes an interesting topic of research in geometric function theory.\vspace{1.2mm}
\par Let $ \mathcal{A} $ denote the set of all analytic functions of the form $ f(z)=\sum_{n=0}^{\infty}a_nz^n $ defined on $ \mathbb{D} $ and we define the class $ \mathcal{B}=\{f\in\mathcal{A} : |f(z)|\leq 1\; \mbox{in}\; \mathbb{D}\} $.  Let us first recall the theorem of Bohr \cite{Bohr-1914} in $ 1914 $, which inspired a lot in the recent years.

\begin{theoA}\cite{Bohr-1914}
	If $ f(z)=\sum_{n=0}^{\infty}a_nz^n\in\mathcal{B} $, then 
	\begin{equation}\label{e-1.1}
		M_f(r):=\sum_{n=0}^{\infty}|a_n|r^n\leq 1 \;\; \mbox{for}\;\; |z|=r\leq\frac{1}{3}.
	\end{equation}
\end{theoA}
Initially, Bohr shows the inequality \eqref{e-1.1} for $|z|\leq1/6$,  but later  M. Riesz, I. Schur and F. Wiener subsequently improved the inequality \eqref{e-1.1} for $|z|\leq1/3$ and showed that the constant $1/3$ is best possible. It is quite natural that the constant $1/3$ and the inequality \eqref{e-1.1} are called respectively, the Bohr radius and Bohr inequality for the class $\mathcal{B}$. Moreover, for:
\begin{equation*}
	f_a(z)=\frac{a-z}{1-az},\;\; a\in [0,1)
\end{equation*} it follows easily that $ M_{f_a}(r)>1 $ if, and only if, $ r>1/(1+2a) $, which shows that $ 1/3 $ is best possible when for $ a\rightarrow 1 $. \vspace{1.2mm}

Let $S_r(f)$ denote the area of the image of the subdisk $\mathbb{D}_r \subseteq \mathbb{D}$ under the mapping $f$. For brevity, and when there is no confusion, we denote this area simply as $S_r$. Recently, Kayumov and Ponnusamy \cite{Kayumov-2017} improved Bohr's inequality \eqref{e-1.1} in various forms. For example, the following inequality was obtained in \cite{Kayumov-CRACAD-2018} (see also \cite{Kayumov-2017}).
\begin{theoB}
	Suppose that $f(z)=\sum_{n=0}^{\infty}a_nz^n$ is analytic in $\mathbb{D}$ and $|f(z)|\leq 1$ in $\mathbb{D}$. Then 
	\begin{align}\label{Eq-1.2}
		|a_0|+\sum_{n=1}^{\infty}|a_n|r^n+\frac{16}{9}\left(\frac{S_r}{\pi}\right)\leq 1\; \mbox{for}\; |z|=r\leq \frac{1}{3}.
	\end{align}
	The numbers $1/3$ and $16/9$ cannot be improved.
\end{theoB}
It is worth recalling that if the first term, $|a_0|$, in equation \eqref{Eq-1.2} is replaced by $|a_0|^2$, then the constants $1/3$ and $16/9$ could be replaced by $1/2$ and $9/8$, respectively (see \cite{Kayumov-CRACAD-2018}). More recently, the authors in \cite{Ismagilov-VINITI-2018} improved Theorem B by replacing the quantity $S_r/\pi$ by $S_r/(\pi-S_r)$. \vspace{1.2mm}

We see that equation \eqref{e-1.1} can be written as
\begin{equation}\label{e-1.2}
	d\left(\sum_{n=0}^{\infty}|a_nz^n|,|a_0|\right)=\sum_{n=1}^{\infty}|a_nz^n|\leq 1-|f(0)|=d(f(0),\partial f(\mathbb{D})),
\end{equation}
where $ d $ is the Euclidean distance and $ \partial f(\mathbb{D}) $ is the boundary of $ f(\mathbb{D}). $ More generally, a class $ \mathcal{F} $ of analytic functions $ f(z)=\sum_{n=0}^{\infty}a_nz^n $ mapping $ \mathbb{D} $ into a domain $ \Omega $ is said to satisfy a Bohr phenomenon if an inequality of type \eqref{e-1.2} holds uniformly in $ |z|\leq \rho_0 $, where $ 0<\rho_0\leq 1 $ for functions in the class $ \mathcal{F} $. Similar definition makes sense for harmonic functions (see \cite{KPSM1018}).\vspace{1.2mm}

\par Let $ \mathcal{H} $ be the class of all complex-valued harmonic functions $ f=h+\bar{g} $ defined on the unit disk $ \mathbb{D} $, where $ h $ and $ g $ are analytic $ \mathbb{D} $ with the normalization $ h(0)=h^{\prime}(0)-1=0 $ and $ g(0)=0 $. Let $ \mathcal{H}_0 $ be defined by $ 	\mathcal{H}_0=\{f=h+\bar{g}\in\mathcal{H} : g^{\prime}(0)=0\}. $ Therefore, each $f=h+\overline{g}\in \mathcal{H}_{0}$ has the following representation 
\begin{equation}\label{e-1.3}
	f(z)=h(z)+\overline{g(z)}=z+\sum_{n=2}^{\infty}a_nz^n+\overline{\sum_{n=2}^{\infty}b_nz^n}.
\end{equation}
Let us recall the Bohr radius for the class $\mathcal{H}_0$ of harmonic mappings.
\begin{defi}
	Let $ f\in\mathcal{H}_0 $ be given by \eqref{e-1.3}. Then the Bohr phenomenon is to find a constant $ R^*\in (0, 1] $ such that the inequality 
	\begin{align*}
		L_f(r):=|z|+\sum_{n=2}^{\infty}\left(|a_n|+|b_n|\right)|z|^n\leq d\left(f(0), \partial\Omega\right) 
	\end{align*} holds for $ |z|=r\leq R^* $, where  $ d\left(f(0), \partial\Omega\right) $ is the Euclidean distance between $ f(0) $ and the boundary of $ \Omega:=f(\mathbb{D}) $. The largest such radius $ R^* $ is called the Bohr radius for the class $ \mathcal{H}_0 $. The radius $R^*$ is said to be best possible if there exists a function $\phi\in  \mathcal{H}_0$ such that $L_{\phi}(r)> d\left(\phi(0), \partial \phi(\mathbb{D})\right)$ for $r>R^*$.
\end{defi}
Stable harmonic mappings were introduced by Hernandez and Martin in \cite{Hernandez-MPCPS-2013}, who also provided definitions for both stable univalent and stable convex harmonic mappings. In \cite{AbdulHadi-BMMSS-2023}, AbdulHadi and Hajj  determined the Bohr radius for the classes of stable univalent harmonic mappings, stable convex harmonic mappings, and stable univalent logharmonic mappings. Recently, in \cite{Liu-Ponn-BSM-2023}, the authors mainly investigated some properties of certain subclasses of stable harmonic mappings on $\mathbb{D}$, including inclusion relations and stability analysis by precise examples, coefficient bounds, growth, covering, and distortion theorems. As applications, they established some Bohr inequalities for these subclasses (stable classes) by means of subordination. For a comprehensive overview of recent developments concerning the Bohr radius problems for various classes of harmonic mappings, readers are referred to the articles \cite{Ahamed-CVEE-2024,Ahamed-Aha-MJM-2024,Aha-Allu-RMJ-2022,Aha-Allu-BMMSS-2022,Ahamed-AMP-2021,Ahamed-CVEE-2021,Aha-Roy-PMS-2025,Liu-Ponn-BSM-2023,Allu-BSM-2021, Evdoridis-IM-2018,Huang-Liu-Ponnusamy-MJM-2021,Liu-Ponnusamy-BMMS-2019} and the references therein.
\begin{defi}
	A (sense-preserving) harmonic mapping $f=h+\overline{g}$ is stable harmonic univalent (SHU) in the unit disk (resp. stable harmonic convex (SHC)) if all the mappings $f_{\mu}=h+\mu\overline{g}$ with $|\mu|=1$ are univalent (resp. convex) in $\mathbb{D}$.
\end{defi}
The harmonic differential operator for functions $f=h+\overline{g}\in\mathcal{S}^0_{\mathcal{H}}$ can be defined as follows:
\begin{defi}
	 If $f=h+\overline{g}\in\mathcal{S}^0_{\mathcal{H}}$, then the harmonic differential operator of $f$ is represented as:
	 \begin{align}\label{Eq-1.4}
	 	Df:=zf_{z}-\overline{z}f_{\overline{z}}
	 \end{align}
	 and 
	 \begin{align}\label{Eq-1.4.1}
	 	\mathscr{D}f:=zf_{z}+\overline{z}f_{\overline{z}}
	 \end{align}
\end{defi}
Furthermore, the linear combination of harmonic mapping and its differential operator is given by 
\begin{align}\label{Eq-1.5}
	F_{\lambda}=(1-\lambda)f+\lambda Df,\; \lambda\in\mathbb{R}.
\end{align}
 \indent In this paper, we establish several results concerning improved Bohr inequalities in terms of distance formulations involving multiple Schwarz functions related to harmonic differential operators $D$ and $\mathscr{D}$ for stable harmonic mappings. The organization of the paper is as follows: In Section \ref{Sec-2}, we present some key lemmas which will play a key role to prove the results of this paper. In Section \ref{Sec-3}, we study improved Bohr inequalities of harmonic differential operators $D$ and $\mathscr{D}$ involving multiple Schwarz functions belong to the class $\mathcal{B}_n$ and obtained Theorems \ref{Th-3.1}, \ref{Th-3.2}, \ref{Th-3.3} and \ref{Th-3.4}. In Section \ref{Sec-4}. we study Bohr phenomenon of the second-order harmonic differential operators $D$ as a convex combinations with the majorent series of harmonic mappings $f$  involving multiple Schwarz functions and establish results Theorems \ref{Th-4.1}, \ref{Th-4.2}, and proved the Theorems \ref{Th-4.3} and \ref{Th-4.4}  for operator $\mathscr{D}^2$ with the majorent series of $f$. Finally, in Section \ref{Sec-5}, we study improved versions of the Bohr inequality for a class of stable harmonic mappings. The improvements are expressed in terms of the quantity ${S_r}/{\pi}$, where $S_r$ is the area of the disk $|z| < r$ under the sense-preserving harmonic map $f_d$ on the unit disk $\mathbb{D}$. Furthermore, we establish results that provide flexibility for a class of harmonic mappings whose associated inequality does not even hold for the class $\mathcal{B}$ of analytic functions, as mentioned in \cite[Remark 2]{Ismagilov-2020-JMAA}, making our results for harmonic mapping  significant.
\section{\bf Auxiliary results}\label{Sec-2}
In this section, we present some key lemmas which will play a crucial role in proving the main results of the paper. These lemmas provide coefficient bounds and growth estimates for the class of stable harmonic mappings and stable convex harmonic mappings, and all were established in \cite{Hernandez-MPCPS-2013}.
\begin{lem}\label{L1}{\cite[Theorem~8.1]{Hernandez-MPCPS-2013}}
Assume that $f = h + \overline{g} \in \mathcal{S}^0_H$ is stable univalent.
Then for all non-negative integers $n$,
\begin{equation}
\big|\,|a_n| - |b_n|\,\big| \leq \max\{|a_n|, |b_n|\} \leq |a_n| + |b_n| \leq n. \tag{2.1}
\end{equation}
These inequalities are sharp. The equalities hold for the analytic Koebe function
\[
k(z) = \frac{z}{(1 - z)^2}.
\]
\end{lem}

\begin{lem}\label{L2}{\cite[Proposition~8.2]{Hernandez-MPCPS-2013}}
Let $f = h + \overline{g} \in \mathcal{S}^0_H$ be a stable harmonic convex mapping.
Then for all non-negative integers $n$,
\begin{equation}
\big|\,|a_n| - |b_n|\,\big| \leq \max\{|a_n|, |b_n|\} \leq |a_n| + |b_n| \leq 1. \tag{2.2}
\end{equation}
These inequalities are sharp. The equalities hold for the function
\[
\ell(z) = \frac{z}{1 - z} \in \mathcal{S}.
\]
\end{lem}

\begin{lem}\label{L3}{\cite[Proposition~8.3]{Hernandez-MPCPS-2013}}
Suppose that $f = h + \overline{g} \in \mathcal{S}^0_H$ is stable harmonic univalent. Then:
\begin{enumerate}[(i)]
\item For all $z \in \mathbb{D}$, we have
\begin{equation}\label{e-2.3}
\frac{|z|}{(1 + |z|)^2} \leq |f(z)| \leq \frac{|z|}{(1 - |z|)^2}. 
\end{equation}
Thus, the disk centered at the origin with radius $\tfrac{1}{4}$ is contained in $f(\mathbb{D})$
for each $f \in \mathcal{S}^0_H$ that is stable harmonic univalent.

\item If $f$ is stable harmonic convex, then
\begin{equation}\label{Eq-2.4}
\frac{|z|}{1 + |z|} \leq |f(z)| \leq \frac{|z|}{1 - |z|}. \tag{2.4}
\end{equation}
Therefore, the disk centered at the origin with radius $\tfrac{1}{2}$ is contained in $f(\mathbb{D})$
for all $f \in \mathcal{S}^0_H$ which are stable convex.
\end{enumerate}
The inequalities in both \textnormal{(i)} and \textnormal{(ii)} are sharp. Equalities hold in \textnormal{(i)} 
if $f$ is a suitable rotation of the Koebe function $k(z)$, and in \textnormal{(ii)} 
for appropriate rotations of the function $\ell(z)$.
\end{lem}
\section{\bf Improved Bohr inequalities of harmonic differential operators involving Schwarz functions}\label{Sec-3}
We introduce the following important subfamily of the class $\mathcal{B}$:
\begin{align*}
	\mathcal{B}_0:=\{\omega\in\mathcal{B} : \omega(0)=0\}=\bigcup_{n=1}^{\infty}\mathcal{B}_n, 
\end{align*} 
where 
\begin{align*}
	\mathcal{B}_n=\bigg\{\omega\in\mathcal{B} : \omega(0)=\cdots=\omega^{(n-1)}(0)=0\; \mbox{and}\; \omega^{(n)}(0)\neq 0\bigg\} \; \mbox{for}\; n\in\mathbb{N}.
\end{align*}
Each function in $\mathcal{B}_0$ is called a Schwarz function. In view of the maximum modulus principle, the family $\mathcal{B}$ consists of non-constant analytic functions $\omega : \mathbb{D}\to\mathbb{D}$ and the constant functions with values in $\overline{\mathbb{D}}$. While some authors exclude constant functions from $\mathcal{B}$, this exclusion does not affect our discussions.\vspace{1.2mm}

In the following research, we will use the series representation of the differential operator. By utilizing equations \eqref{e-1.3} and \eqref{Eq-1.4}, we derive the expansion given below
\begin{align}	\label{m1}
	Df = z + \sum_{n=2}^{\infty} n a_n z^n - \ol{\sum_{n=2}^{\infty} n b_n z^n}.
\end{align}
and
\begin{align}\label{m2}
\mathscr{D}f(z) = z + \sum_{n=2}^{\infty} n a_n z^n + \overline{\sum_{n=2}^{\infty} n b_n z^n}
\end{align}
We now state our first result for $f = h + \overline{g} \in \mathcal{S}^0_H$, where $f$ is a stable univalent harmonic mapping on the unit disk $\mathbb{D}$. This result presents an improved Bohr inequality involving multiple Schwarz functions, specifically $\omega_m\in \mathcal{B}_m$ and $\omega_p\in \mathcal{B}_p$.
\begin{theo}\label{Th-3.1}
Let $f = h + \overline{g} \in \mathcal{S}^0_H$ be a stable univalent harmonic mapping on the unit disk $\mathbb{D}$, and let $\omega_m\in \mathcal{B}_m$ and $ \omega_p\in \mathcal{B}_p$. Then 
\begin{align*}
	&|Df(\omega_m(z))| + M_f(|\omega_p(z)|) \leq d\big(f(0), \partial f(\mathbb{D})\big)
\end{align*}
and 
\begin{align*}
	|\mathscr{D}f(\omega_m(z))| + M_f(|\omega_p(z)|) \leq d\big(f(0), \partial f(\mathbb{D})\big)
\end{align*}
both hold for $|z|=r\leq r_{1,m,p}$, where $r_{1,m,p}$ is the unique root in $(0,1)$ of the equation
\begin{align}\label{rmp}
	r^m+r^p+\frac{r^{2m}(4-3r^{m}+r^{2m})}{(1-r^m)^3}+\frac{2r^{2p}-r^{3p}}{(1-r^p)^2}-\frac{1}{4}=0.
\end{align}
The radius $r_{1,m,p}$ is the best possible.
\end{theo}
\begin{rem}
	In particular, taking $\omega_m(z)=z$ and $\omega_p(z) = z$ for $m=1=p$, Theorem \ref{Th-3.1} simplifies to \cite[Theorem 3.1]{Xu-An-Liu-JA-2025}. This confirms that our result provides a significant generalization and improvement.
\end{rem}
\begin{table}[H]
\centering
\begin{minipage}[b]{0.45\textwidth}
\centering
\begin{tabular}{|c|c|c|}
\hline
\textbf{m} & \textbf{p} & \textbf{$r_{1,m,p}$} \\ 
\hline
1 & 1 & 0.093200 \\ 
\hline
2 & 1 & 0.157800 \\ 
\hline
2 & 2 & 0.305300 \\ 
\hline
1 & 2 & 0.133100 \\ 
\hline
\end{tabular}
\caption{$r_{1,m,p}$ is the smallest root of equation (\ref{rmp}) in $(0,1)$.}
\label{tab3.1}
\end{minipage}%
\hfill
\begin{minipage}[b]{0.55\textwidth}
\centering
\begin{figure}[H]
\centering
\includegraphics[width=\textwidth]{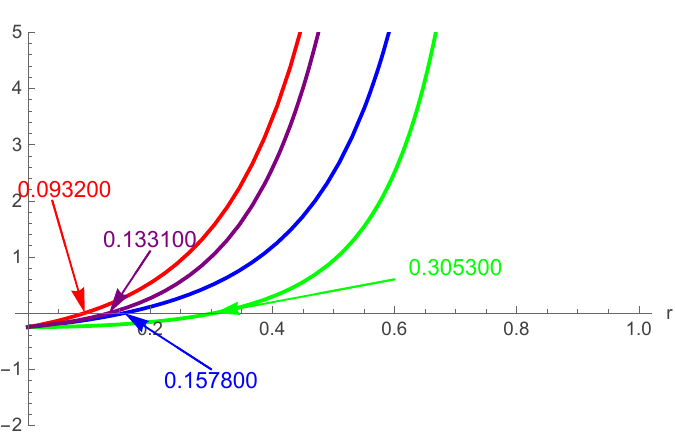}
\caption{The graphs exhibit the locations of the roots $r_{1,m,p}$ in $(0,1)$ for different values of $m,p$.}
\label{f31}
\end{figure}
\end{minipage}
\end{table}

\begin{proof}[\bf Proof of Theorem \ref{Th-3.1}]
If $f$ is a stable univalent harmonic mapping, then from (2.3) we have
\[
|f(z)| \ge \frac{|z|}{(1 + |z|)^2}.
\]
It follows from equation \eqref{e-2.3} of Lemma \ref{L3} that the Euclidean distance between $f(0)$ and the boundary of the image domain $f(\mathbb{D})$ satisfies the inequality
\begin{align}\label{t1.1}
	d\big(f(0), \partial f(\mathbb{D})\big)
	= \lim_{|z| \to 1} \inf |f(z) - f(0)|
	\ge \lim_{|z| \to 1} \inf \frac{|z|}{(1 + |z|)^2}
	= \frac{1}{4}.
\end{align}
For $m \ge 1$, the condition $|\omega_m(z)| \leq |z|^m$ holds. By setting $|z|=r$ and applying this inequality (from (\ref{m1})) in conjunction with Lemma \ref{L1}, a straightforward computation shows that
\begin{align*}
|Df(\omega_m(z))| + M_f(|\omega_p(z)|)
&\le |\omega_m(z)|
   + \sum_{n=2}^{\infty} n |a_n| |\omega_m(z)|^n
   + \sum_{n=2}^{\infty} n |b_n| |\omega_m(z)|^n  \\
&\quad + |\omega_p(z)|
   + \sum_{n=2}^{\infty} (|a_n| + |b_n|) |\omega_p(z)|^n \\
&\le |z|^m + \sum_{n=2}^{\infty} n (|a_n| + |b_n|) |z|^{mn}
     + |z|^p + \sum_{n=2}^{\infty} (|a_n| + |b_n|) |z|^{pn} \\
&\le r^m + \sum_{n=2}^{\infty} n^2 r^{mn}
      + r^p + \sum_{n=2}^{\infty} n r^{np} \\
&=\frac{1}{4}+G_1(r), 
\end{align*}
where 
\begin{align*}
	G_1(r):= r^m + \frac{r^{2m}(4 - 3r^{m} + r^{2m})}{(1 - r^m)^3}+ r^p + \frac{2r^{2p} - r^{3p}}{(1 - r^p)^2}-\frac{1}{4}.
\end{align*}
It follows directly from the definition that $G_1(r)$ is differentiable on $[0,1).$ Differentiating $G_1(r)$ with respect to $r$, we obtain
\begin{align*}
	G_1'(r)&=mr^{m-1}+\frac{(8r^m-9r^{2m}+4r^{3m})mr^{m-1}}{(1 - r^m)^3}+\frac{3r^{2m}(4 - 3r^{m} + r^{2m})}{(1 - r^m)^4}mr^{m-1}\\
	&\quad+pr^{p-1}+\frac{p(4r^{2p-1} - 3r^{3p-1})}{(1 - r^p)^2}+\frac{2p(2r^{2p} - r^{3p})r^{p-1}}{(1 - r^p)^3}\\
	&= \frac{m(1+5r^m-6r^{2m}+3r^{3m})r^{m-1}}{(1 - r^m)^3}+\frac{3r^{2m}(4 - 3r^{m} + r^{2m})}{(1 - r^m)^4}mr^{m-1}\\
	&\quad +pr^{p-1}+\frac{4p(r^{2p-1} - 3r^{3p-1})}{(1 - r^p)^2}+\frac{2p(2r^{2p} - r^{3p})r^{p-1}}{(1 - r^p)^3}.
\end{align*}
Since all terms of $G_1'(r)$ are positive, $G_1(r)$ is strictly increasing on $[0,1)$. Thus we have $G_1(r)\geq G_1(r_{m,p})$ for $r\geq r_{m,p}$ and this implies the inequality
\begin{align}\label{Eq-3.3}
	r^m &+ \frac{r^{2m}(4 - 3r^{m} + r^{2m})}{(1 - r^m)^3}+ r^p + \frac{2r^{2p} - r^{3p}}{(1 - r^p)^2}\\&\geq r_{m,p}^m + \frac{r_{m,p}^{2m}(4 - 3r_{m,p}^{m} + r_{m,p}^{2m})}{(1 - r_{m,p}^m)^3}+ r_{m,p}^p + \frac{2r_{m,p}^{2p} - r_{m,p}^{3p}}{(1 - r_{m,p}^p)^2}.\nonumber
\end{align}
A simple calculation shows that 
\begin{align*}
	G_1(0)=-\frac{1}{4}<0,\;\;\text{and}\;\;\lim_{r\rightarrow 1^-}G_1(r)=+\infty.
\end{align*}
Consequently, $G_1(r)$ has a unique root in the interval $(0,1)$. Denoting this unique root by $r_{m,p}$, we obtain the desired inequality
\begin{align*}
	|Df(\omega_m(z))| + M_f(|\omega_p(z)|)
	\le \frac{1}{4}
	\le d\big(f(0), \partial f(\mathbb{D})\big),
\end{align*}
for $|z| \le r_{m,p}$.\vspace{1.2mm}

Next part of the proof is to show the radius $r_{m,p}$ is best possible. Henceforth, we consider the Koebe function
\begin{align}	\label{f1}
	k(z) = \frac{z}{(1 - z)^2} = z + \sum_{n=2}^{\infty} n z^n,
\end{align}
which is harmonic and univalent in $\mathbb{D}$. It follows from \eqref{e-2.3} of Lemma \ref{L3} that 
\begin{equation}\label{Eq-3.4}
	d\big(k(0), \partial k(\mathbb{D})\big)=\frac{1}{4}.
\end{equation}  
Considering $\omega_m(z) = z^m$ ($m\geq 1$) and the case where $z = r>r_{m,p}$ with $f=k$, it follows from equations \eqref{Eq-3.3} and \eqref{Eq-3.4} that
\begin{align*}
	|Df(\omega_m(z))| + M_f(|\omega_p(z)|)&=|Dk(r^m)| + M_{k}(r^p)\\
	&= r^m + \sum_{n=2}^{\infty} n^2 r^{mn}
	+ r^p + \sum_{n=2}^{\infty} n r^{np}\\
	&= r^m + \frac{r^{2m}(4 - 3r^{m} + r^{2m})}{(1 - r^m)^3}+ r^p + \frac{2r^{2p} - r^{3p}}{(1 - r^p)^2}\\&>r_{m,p}^m + \frac{r_{m,p}^{2m}(4 - 3r_{m,p}^{m} + r_{m,p}^{2m})}{(1 - r_{m,p}^m)^3}+ r_{m,p}^p + \frac{2r_{m,p}^{2p} - r_{m,p}^{3p}}{(1 - r_{m,p}^p)^2}\\&=\frac{1}{4}\\&=d\big(k(0), \partial k(\mathbb{D})\big)
\end{align*}
which shows that the radius $r_{m,p}$ is the best possible.\vspace{1.2mm}

Similarly, for the operator $\mathscr{D}$, we obtain the inequality
\begin{align*}
	|\mathscr{D}f(\omega_m(z))| + M_f(|\omega_p(z)|)
	&\le |\omega_m(z)|
	+ \sum_{n=2}^{\infty} n |a_n| |\omega_m(z)|^n
	+ \sum_{n=2}^{\infty} n |b_n| |\omega_m(z)|^n  \\
	&\quad + |\omega_p(z)|
	+ \sum_{n=2}^{\infty} (|a_n| + |b_n|) |\omega_p(z)|^n \\
	&\le \frac{1}{4}+G_1(r).
\end{align*}
holds. The remaining part of the proof is analogous to the argument used for the inequality involving the operator $D$, and thus, we omit the details. This completes the proof.
\end{proof}
We now state our second result for $f = h + \overline{g} \in \mathcal{S}^0_H$, where $f$ is a stable convex harmonic mapping on the unit disk $\mathbb{D}$. This result presents an improved Bohr inequality involving multiple Schwarz functions $\omega_m\in \mathcal{B}_m$ and $\omega_p\in \mathcal{B}_p$.
\begin{theo}\label{Th-3.2}
Let $f = h + \overline{g} \in \mathcal{S}^0_H$ be a stable convex harmonic mapping on the unit disk $\mathbb{D}$ and $\omega_m\in \mathcal{B}_m$ and $\omega_p\in \mathcal{B}_p$. Then 
\begin{align*}
&|Df(\omega_m(z))| + M_f(|\omega_p(z)|) \leq d\big(f(0), \partial f(\mathbb{D})\big)
\end{align*}
and 
\begin{align*}
	|\mathscr{D}f(\omega_m(z))| + M_f(|\omega_p(z)|) \leq d\big(f(0), \partial f(\mathbb{D})\big)
\end{align*}
both hold for $|z|=r \leq r_{2,m,p}$, where $r_{2,m,p}$ is the unique root in $(0,1)$ of the equation
\bea\label{rmp2}
r^m+r^p+\frac{2r^{2m} - r^{3m}}{(1 - r^m)^2}+\frac{r^{2p}}{1-r^p}-\frac{1}{2}=0.
\eea
The radius $r_{2,m,p}$ is the best possible.
\end{theo}
\begin{rem}
	In particular, taking $\omega_m(z)=z$ and $\omega_p(z) = z$ for $m=1=p$, Theorem \ref{Th-3.2} simplifies to \cite[Corollary 3.2]{Xu-An-Liu-JA-2025}. This confirms that our result provides a significant generalization and improvement.
\end{rem}
\begin{table}[H]
\centering
\begin{minipage}[b]{0.45\textwidth}
\centering
\begin{tabular}{|c|c|c|}
\hline
\textbf{m} & \textbf{p} & \textbf{$r_{2,m,p}$} \\ 
\hline
1 & 1 & 0.183500 \\ 
\hline
2 & 1 & 0.386900 \\ 
\hline
2 & 2 & 0.428400 \\ 
\hline
1 & 2 & 0.246800 \\ 
\hline
\end{tabular}
\caption{$r_{2,m,p}$ is the smallest  root of equation (\ref{rmp2}) in $(0,1)$.}
\label{tab3.2}
\end{minipage}%
\hfill
\begin{minipage}[b]{0.55\textwidth}
\centering
\begin{figure}[H]
\centering
\includegraphics[width=\textwidth]{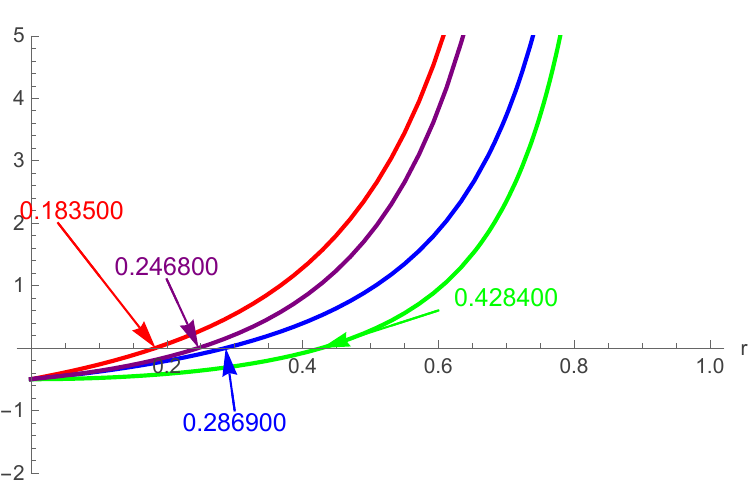}
\caption{The graphs exhibit the locations of the roots $r_{2,m,p}$ in $(0,1)$ for different values of $m,p$.}
\label{f32}
\end{figure}
\end{minipage}
\end{table}

\begin{proof}[\bf Proof of Theorem \ref{Th-3.2}]
 Let $f = h + \overline{g} \in \mathcal{S}^0_H$ be a stable convex harmonic mapping on the unit disk $\mathbb{D}$. Then it follows from \eqref{Eq-2.4} of Lemma \ref{L3} that 
\[
|f(z)| \ge \frac{|z|}{1 + |z|}.
\]
Thus, the Euclidean distance between $f(0)$ and $\partial f(\mathbb{D})$ satisfies the inequality
\begin{equation}\label{t2.1}
d\big(f(0), \partial f(\mathbb{D})\big)
  = \lim_{|z| \to 1} \inf |f(z) - f(0)|
  \ge \lim_{|z| \to 1} \inf \frac{|z|}{1 + |z|}
  = \frac{1}{2}.
\end{equation}
Thus, for $|z| = r$, combining (\ref{m1}) with Lemma~\ref{L2}, a simple computation shows that 
\begin{align*}
|Df(\omega_m(z))| + M_f(|\omega_p(z)|)
&\le |\omega_m(z)|
   + \sum_{n=2}^{\infty} n |a_n| |\omega_m(z)|^n
   + \sum_{n=2}^{\infty} n |b_n| |\omega_m(z)|^n  \\
&\quad + |\omega_p(z)|
   + \sum_{n=2}^{\infty} (|a_n| + |b_n|) |\omega_p(z)|^n \\
&\le |z|^m + \sum_{n=2}^{\infty} n (|a_n| + |b_n|) |z|^{mn}
     + |z|^p + \sum_{n=2}^{\infty} (|a_n| + |b_n|) |z|^{pn} \\
&\le r^m + \sum_{n=2}^{\infty} n r^{mn}
      + r^p + \sum_{n=2}^{\infty} r^{np} \\
&=\frac{1}{2}+G_2(r), 
\end{align*}
where 
\begin{align*}
	G_2(r)=r^m+r^p+\frac{2r^{2m} - r^{3m}}{(1 - r^m)^2}+\frac{r^{2p}}{1-r^p}-\frac{1}{2}.
\end{align*}
It is evident that $G_2(r)$ is continuous on $[0,1)$. Note that 
\begin{align*}
	G_2(0)=-\frac{1}{2}<0\;\;\text{and}\;\;\lim_{r\rightarrow 1^{-}}G_2(r)=+\infty.
\end{align*}
Differentiating $G_2(r)$, we see that 
\begin{align*}
	G_2'(r)&=mr^{m-1}+pr^{p-1} +\frac{m(4r^{2m} - 3r^{3m})}{(1 - r^m)^2}+\frac{2(2r^{2m} - r^{3m})}{(1 - r^m)^3}\\&\quad+\frac{2pr^{2p-1}}{1-r^p}+\frac{r^{2p}}{(1-r^p)^2}> 0.
\end{align*}
This shows that $G_2(r)$ is increasing function on $[0,1)$.  Consequently, $G_1(r)$ admits a unique root in $(0,1)$. Let $r_{2,m,p}$ is the unique root of $G_2(r)$ in $(0,1)$. then for $r\geq r_{2,m,p}$,  we have 
\begin{align}\label{Eq-3.7}
	r^m&+r^p+\frac{2r^{2m} - r^{3m}}{(1 - r^m)^2}+\frac{r^{2p}}{1-r^p}-\frac{1}{2}\\&\geq r_{2,m,p}^m+r_{2,m,p}^p+\frac{2r_{2,m,p}^{2m} - r_{2,m,p}^{3m}}{(1 - r_{2,m,p}^m)^2}+\frac{r_{2,m,p}^{2p}}{1-r_{2,m,p}^p}.\nonumber
\end{align}

Thus, we conclude that the inequality
\[
|Df(\omega_m(z))| + M_f(|\omega_p(z)|)
   \le \frac{1}{2}
   \le d\big(f(0), \partial f(\mathbb{D})\big),
\]
holds for $|z| \le r_{2,m,p}$.\vspace{1.2mm}

To establish the radius $r_{2,m,p}$ is best possible, we consider the analytic convex function
\begin{align}\label{f2}
	f_2(z) = \frac{z}{1 - z} = z + \sum_{n=2}^{\infty} z^n,
\end{align}
which is harmonic and univalent in $\mathbb{D}$. It follows from \eqref{e-2.3} of Lemma \ref{L3} that 
\begin{equation}\label{Eq-3.8}
	d\big(f_2(0), \partial f_2(\mathbb{D})\big)=\frac{1}{2}.
\end{equation}  
 For $\omega_m(z) = z^m$ with $m \ge 1$, $z = r>r_{2,m,p}$ and $f=f_2$, a simple computation using \eqref{Eq-3.7} and \eqref{Eq-3.8} yields
\begin{align*}
	|Df(\omega_m(z))| + M_f(|\omega_p(z)|)&=|D{f_1}(r^m)| + M_{f_1}(r^p)\\
	&= r^m + \sum_{n=2}^{\infty} n r^{mn}
	+ r^p + \sum_{n=2}^{\infty}  r^{np}\\
	&= r^m+r^p+\frac{2r^{2m} - r^{3m}}{(1 - r^m)^2}+\frac{r^{2p}}{1-r^p}\\&>r_{2,m,p}^m+r_{2,m,p}^p+\frac{2r_{2,m,p}^{2m} - r_{2,m,p}^{3m}}{(1 - r_{2,m,p}^m)^2}+\frac{r_{2,m,p}^{2p}}{1-r_{2,m,p}^p}\\&=\frac{1}{2}\\&=d\big(f_2(0), \partial f_2(\mathbb{D})\big).
\end{align*}
This establishes $r_{2,m,p}$ as the best possible radius.

Similarly, for the operator $\mathscr{D}$, we obtain the inequality
\begin{align*}
	|\mathscr{D}f(\omega_m(z))| + M_f(|\omega_p(z)|)&\le \frac{1}{2}+G_2(r).
\end{align*}
holds. The remaining part of the proof is analogous to the argument used for the inequality involving the operator $D$, and thus, we omit the details.  this completes the proof.
\end{proof}
In the following result, we investigate an improved Bohr-type inequality for stable univalent harmonic mappings associated with differential operators. Considering $ f = h + \overline{g} \in \mathcal{S}_H^0 $, we establish sharp bounds involving the functions $ Df $ and $\mathscr{D}f $ in terms of multiple Schwarz functions $\omega_m \in \mathcal{B}_m$ and $\omega_q \in \mathcal{B}_q$.
\begin{theo}\label{Th-3.3}
Let $f = h + \overline{g} \in \mathcal{S}^0_H$ be a stable univalent harmonic mapping on the unit disk $\mathbb{D}$ and $\omega_m\in \mathcal{B}_m$ and $\omega_q\in \mathcal{B}_q$. Then for $s, p\in\mathbb{N}$, we have
\begin{align*}
	|Df(\omega_m(z))|^s +|f(\omega_q(z))|^p\leq 1\;\;\text{and}\;\;|\mathscr{D}f(\omega_m(z))|^s +|f(\omega_q(z))|^p\leq 1
\end{align*}
for $|z|=r \leq r_{3,s,m,p,q}$, where $r_{3,s,m,p,q}$ is the unique root in $(0,1)$ of the equation
\begin{align}\label{rmp3}
	\frac{r^{sm}(1+r^{m})^s}{(1-r^m)^{3s}}+\left(\frac{r^q}{(1-r^q)^2}\right)^p-1=0.
\end{align}
The radius $r_{3,s,m,p,q}$ is the best possible.
\end{theo}
\begin{rem}
By setting $\omega_m(z) = z$ and $\omega_p(z) = z$ with $m = p = 1$, $s=2$, Theorem \ref{Th-3.3} reduces to \cite[Theorem 3.3]{Xu-An-Liu-JA-2025}. Hence, our theorem extends and improves the earlier result.
\end{rem}
\begin{table}[H]
\centering
\begin{minipage}[b]{0.45\textwidth}
\centering
\begin{tabular}{|c|c|c|c|c|}
\hline
\textbf{s} &\textbf{m} & \textbf{p} &\textbf{q} &\textbf{$r_{3,s,m,p,q}$} \\ 
\hline
2& 1 & 1 & 1&0.250500 \\ 
\hline
2 & 2 & 3& 1& 0.378000 \\ 
\hline
3 & 2 & 5& 5 & 0.533600 \\ 
\hline
2 & 1 & 7& 2 & 0.284800 \\ 
\hline
\end{tabular}
\caption{$r_{3,s,m,p,q}$ is the smallest root of equation (\ref{rmp3}) in $(0,1)$.}
\label{tab3.3}
\end{minipage}%
\hfill
\begin{minipage}[b]{0.55\textwidth}
\centering
\begin{figure}[H]
\centering
\includegraphics[width=\textwidth]{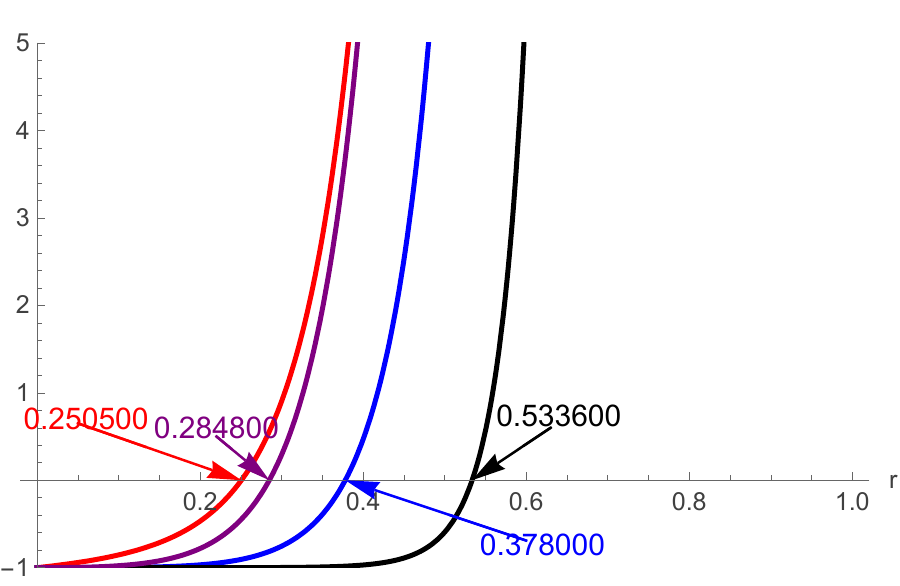}
\caption{The graphs exhibit the locations of the roots $r_{3,s,m,p,q}$ in $(0,1)$ for different values of $s,m,p,q$.}
\label{f33}
\end{figure}
\end{minipage}
\end{table}

\begin{proof}[\bf Proof of Theorem \ref{Th-3.3}]  By Lemma \ref{L1}, Lemma \ref{L3} (i) and (\ref{m1}), we have
\begin{align}	\label{t3.1} |Df(\omega_m(z))|^s +|f(\omega_q(z))|^p &\leq \left(|\omega_m(z)|+ \sum_{n=2}^{\infty} n |a_n| |\omega_m(z)|^n+ \sum_{n=2}^{\infty} n |b_n| |\omega_m(z)|^n\right)^s\nonumber\\
	&\quad+\left(\frac{|\omega_q(z)|}{(1-|\omega_q(z)|)^2}\right)^p.
\end{align}
As $\omega_q\in \mathcal{B}_q$ in view of the classical Schwarz Lemma, we have $|\omega_q(z)|\leq |z|^q$; furthermore, we see that 
\begin{align}
	\label{t3.2} \frac{|\omega_q(z)|}{(1-|\omega_q(z)|)^2}&=|\omega_q(z)|(1+2|\omega_q(z)|+3|\omega_q(z)|^2+\cdots)\nonumber\\
	&\leq |z|^q(1+2|z|^q+3|z|^{2q}+\cdots)\nonumber\\
	&=\frac{|z|^q}{(1-|z|^q)^2}.
\end{align}
For $|z|=r$, it follows from (\ref{t3.3}), (\ref{t3.2}), and the Schwarz Lemma that 
\begin{align}
	\label{t3.3} |Df(\omega_m(z))|^s +|f(\omega_q(z))|^p &\leq \left(|z|^m+ \sum_{n=2}^{\infty} n (|a_n| +|b_n|) |z|^{mn}\right)^s+\left(\frac{|z|^q}{1-|z|^q}\right)^p\nonumber\\
	&= \left(r^m+ \sum_{n=2}^{\infty} n (|a_n|+|b_n|) r^{mn}\right)^s+\left(\frac{r^q}{(1-r^q)^2}\right)^p\nonumber\\
	&\leq  \left(r^m+ \sum_{n=2}^{\infty} n^2 r^{mn}\right)^s+\left(\frac{r^q}{(1-r^q)^2}\right)^p\nonumber\\
	&=\left(r^m+\frac{r^{2m}(4-3r^m+r^{2m})}{(1-r^m)^3}\right)^s+\left(\frac{r^q}{(1-r^q)^2}\right)^p\nonumber\\
	&=\frac{r^{sm}(1+r^{m})^s}{(1-r^m)^{3s}}+\left(\frac{r^q}{(1-r^q)^2}\right)^p\nonumber\\
	&=1+G_3(r),
\end{align}
where 
\begin{align*}
	G_3(r)=\frac{r^{sm}(1+r^{m})^s}{(1-r^m)^{3s}}+\left(\frac{r^q}{(1-r^q)^2}\right)^p-1.
\end{align*}
It is evident that $G_3(r)$ is continuous on the interval $[0,1)$. Moreover,  
\begin{align*}
	G_3(0) = -1 < 0, 
	\quad \text{and} \quad 
	\lim_{r \to 1^{-}} G_3(r) = +\infty.
\end{align*}
Differentiating $G_3(r)$ with respect to $r$, we obtain
\[
G_3'(r)
=\frac{s m\, r^{s m - 1}(1 + r^m)^{s - 1}(1 + 2r^m)}{(1 - r^m)^{3s + 1}}+\frac{p q\, r^{p q - 1} (1 + r^q)}{(1 - r^q)^{2p + 1}}>0.
\]
Hence, $G_3(r)$ is a monotonically increasing function on $[0,1)$.  
Consequently, $G_3(r)$ admits a unique zero in the interval $(0,1)$.  
Let $r_{2,s,m,p,q}$ denote this unique root of $G_3(r)$. Then, for $r\geq r_{2,s,m,p,q}$, we see that 
\begin{align}\label{Eq-3.13}
	\frac{r^{sm}(1+r^{m})^s}{(1-r^m)^{3s}}+\left(\frac{r^q}{(1-r^q)^2}\right)^p\geq \frac{r_{2,s,m,p,q}^{sm}(1+r_{2,s,m,p,q}^{m})^s}{(1-r_{2,s,m,p,q}^m)^{3s}}+\left(\frac{r_{2,s,m,p,q}^q}{(1-r_{2,s,m,p,q}^q)^2}\right)^p.
\end{align} 
Thus, the desired inequality holds
\begin{align*}
	|Df(\omega_m(z))|^s + |f(\omega_q(z))|^p\leq 1,
\end{align*}
for $|z| \le r_{3,s,m,p,q}$.\vspace{1.2mm}

To show the radius $r_{3,s,m,p,q}$ is best possible, we consider the analytic function
\begin{align*}
	f_3(z) = \frac{z}{(1 - z)^2} = z + \sum_{n=2}^{\infty} nz^n,
\end{align*}
which is harmonic and univalent in $\mathbb{D}$. \vspace{2mm}
 
For $\omega_m(z) = z^m$ (where $m \ge 1$), $z = r>r_{3,s,m,p,q}$, and $f=f_3$, a straightforward computation using \eqref{Eq-3.13} shows that
\begin{align*}
	|Df(\omega_m(z))|^s + |f(\omega_q(z))|^p&=|D_{f_3}(r^m)|^s +|f(r^q)|^p\\&=
	\frac{r^{sm}(1+r^{m})^s}{(1-r^m)^{3s}}+\left(\frac{r^q}{1-r^q}\right)^p\\&>\frac{r_{3,s,m,p,q}^{sm}(1+r_{3,s,m,p,q}^{m})^{s}}{(1-r_{3,s,m,p,q}^m)^{3s}}+\left(\frac{r_{3,s,m,p,q}^q}{(1-r_{3,s,m,p,q}^q)^2}\right)^p
	\\&= 1.
\end{align*}
This shows that the radius $r_{3,s,m,p,q}$ is the best possible.\vspace{1.2mm}

Similarly, for the operator $\mathscr{D}$, we can easily shown that  the inequality
\begin{align*}
|\mathscr{D}f(\omega_m(z))|^s + |f(\omega_q(z))|^p&\le \frac{1}{2}+G_3(r).
\end{align*}
holds. The remaining part of the proof is analogous to the argument used for the inequality involving the operator $D$, and thus, we omit the details. This completes the proof.
\end{proof}
For $f = h + \overline{g} \in \mathcal{S}^0_H$, a stable convex harmonic mapping in the unit disk $\mathbb{D}$, we obtain the improved Bohr inequality considering integral powers involving multiple Schwarz functions $\omega_m \in \mathcal{B}_m$ and $\omega_q \in \mathcal{B}_q$.
\begin{theo}\label{Th-3.4}
Let $f = h + \overline{g} \in \mathcal{S}^0_H$ be a stable convex harmonic mapping in the unit disk $\mathbb{D}$, and let $\omega_m \in \mathcal{B}_m$ and $\omega_q \in \mathcal{B}_q$. Then, for $s, p \in\mathbb{N}$, we have
\begin{align*}|Df(\omega_m(z))|^s +|f(\omega_q(z))|^p\leq 1\;\;\text{and}\;\;|\mathscr{D}f(\omega_m(z))|^s +|f(\omega_q(z))|^p\leq 1
\end{align*}
for $|z|=r \leq r_{4,s,m,p,q}$, where $r_{4,s,m,p,q}$ is the unique root in $(0,1)$ of the equation
\begin{equation}\label{rmp4}
\frac{r^{sm}}{(1 - r^m)^{2s}}
+ \left(\frac{r^q}{1 - r^q}\right)^p - 1 = 0.
\end{equation}
The radius $r_{4,s,m,p,q}$ is the best possible.
\end{theo}
\begin{rem}
For the special case $\omega_m(z) = z$ and $\omega_p(z) = z$ with $m = p = 1$, $s=2$, Theorem \ref{Th-3.4} reduces to \cite[Corollary 3.4]{Xu-An-Liu-JA-2025}. Hence, our theorem improves the earlier result.
\end{rem}
\begin{table}[H]
\centering
\begin{minipage}[b]{0.45\textwidth}
\centering
\begin{tabular}{|c|c|c|c|c|}
\hline
\textbf{s} &\textbf{m} & \textbf{p} &\textbf{q} &\textbf{$r_{4,s,m,p,q}$} \\ 
\hline
2& 1 & 1 & 1&0.326200 \\ 
\hline
2 & 2 & 3& 1& 0.485200 \\ 
\hline
3 & 2 & 5& 5 & 0.618300 \\ 
\hline
2 & 1 & 7& 2 & 0.382000 \\ 
\hline
\end{tabular}
\caption{$r_{4,s,m,p,q}$ is the smallest root of equation (\ref{rmp4}) in $(0,1)$.}
\label{tab3.4}
\end{minipage}%
\hfill
\begin{minipage}[b]{0.55\textwidth}
\centering
\begin{figure}[H]
\centering
\includegraphics[width=\textwidth]{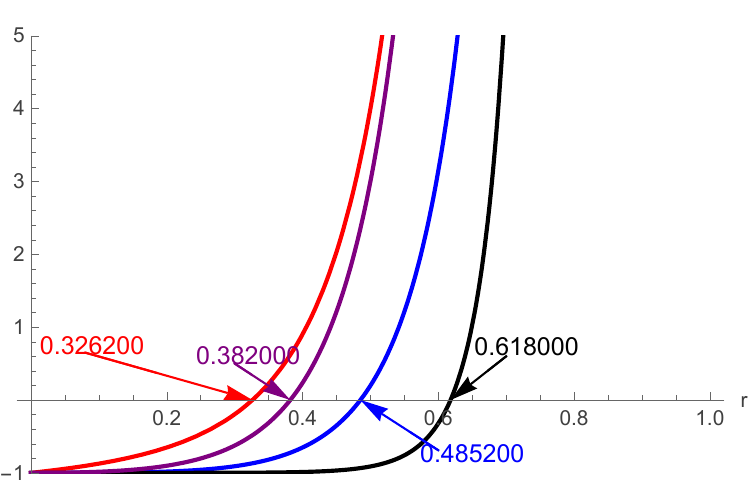}
\caption{The graphs exhibit the locations of the roots $r_{4,s,m,p,q}$ in $(0,1)$ for different values of $s,m,p,q$.}
\label{f34}
\end{figure}
\end{minipage}
\end{table}

\begin{proof}
In view of (\ref{m1}), and by Lemmas \ref{L1} and \ref{L3}(i), we obtain
\begin{align}\label{t4.1}
|Df(\omega_m(z))|^s + |f(\omega_q(z))|^p
&\le \left(|\omega_m(z)| + \sum_{n=2}^{\infty} n |a_n| |\omega_m(z)|^n
   + \sum_{n=2}^{\infty} n |b_n| |\omega_m(z)|^n \right)^s \nonumber\\
&\quad + \left(\frac{|\omega_q(z)|}{1 - |\omega_q(z)|}\right)^p.
\end{align}
As $\omega_q\in \mathcal{B}_q$ in view of the classical Schwarz Lemma, we have $|\omega_q(z)|\leq |z|^q$; furthermore, we see that 
\begin{align}
	\label{ttt.2} \frac{|\omega_q(z)|}{1-|\omega_q(z)|}&=|\omega_q(z)|(1+|\omega_q(z)|+|\omega_q(z)|^2+\cdots)\nonumber\\
	&\leq |z|^q(1+|z|^q+|z|^{2q}+\cdots)\nonumber\\
	&=\frac{|z|^q}{1-|z|^q}.
\end{align}
For $|z| = r$, using (\ref{t4.1}), (\ref{ttt.2}), and the Schwarz lemma, we get
\begin{align}\label{t4.3}
|Df(\omega_m(z))|^s + |f(\omega_q(z))|^p
&\le \left(r^m + \sum_{n=2}^{\infty} n(|a_n| + |b_n|) r^{mn}\right)^s
     + \left(\frac{r^q}{1 - r^q}\right)^p \nonumber\\
&\le \left(r^m + \sum_{n=2}^{\infty} n r^{mn}\right)^s
     + \left(\frac{r^q}{1 - r^q}\right)^p \nonumber\\
&= \left(r^m + \frac{r^{2m}(2 -  r^{m})}{(1 - r^m)^2}\right)^s
   + \left(\frac{r^q}{1 - r^q}\right)^p \nonumber\\
&= \frac{r^{sm}}{(1 - r^m)^{2s}}
   + \left(\frac{r^q}{1 - r^q}\right)^p \nonumber\\
&= 1 + G_4(r),
\end{align}
where
\[
G_4(r) = \frac{r^{sm}}{(1 - r^m)^{2s}} + \left(\frac{r^q}{1 - r^q}\right)^p - 1.
\]

It is evident that $G_3(r)$ is continuous on $[0,1)$.  
Moreover,
\[
G_3(0) = -1 < 0\; \mbox{and}\;
\lim_{r \to 1^{-}} G_3(r) = +\infty.
\]
Differentiating $G_3(r)$ with respect to $r$ gives
\begin{align*}
	G_4'(r)
	= \frac{r^{sm}}{(1 - r^m)^{2s}} \,
	sm\!\left(\frac{1}{r} + \frac{2r^{m-1}}{1 - r^m}\right)
	+ pq\,r^{q-1}
	\left(\frac{r^q}{1 - r^q}\right)^{p-1}
	\frac{1}{(1 - r^q)^2}\geq 0
\end{align*}
Hence, $G_4(r)$ is a monotonically increasing function on $[0,1)$.
Therefore, $G_4(r)$ admits a unique zero in $(0,1)$.
Let $r_{4,s,m,p,q}$ denote this unique root in $(0,1)$. Consequently, the inequality
\begin{align*}
	|Df(\omega_m(z))|^s + |f(\omega_q(z))|^p \le 1
\end{align*}
holds for all $|z| \le r_{4,s,m,p,q}$.\vspace{1.2mm}

To show the radius $r_{4,s,m,p,q}$ is best possible, we consider the analytic function 
\begin{align*}
	f_4(z) = \frac{z}{1 - z}
	= z + \sum_{n=2}^{\infty} z^n,
\end{align*}
which is convex harmonic and univalent in $\mathbb{D}$.  \vspace{1.2mm}

For $\omega_m(z) = z^m$ (where $m \ge 1$), $z = r>r_{4,s,m,p,q}$, and $f=f_4$, an easy computation leads to
\begin{align*}
	|Df(\omega_m(z))|^s + |f(\omega_q(z))|^p&=|D{f_3}(r^m)|^s + |f_3(r^q)|^p\\
	&= \frac{r^{sm}}{(1 - r^m)^{2s}}
	+ \left(\frac{r^q}{1 - r^q}\right)^p\\&>\frac{r_{4,s,m,p,q}^{sm}}
	{(1 - r_{4,s,m,p,q}^m)^{2s}}
	+ \left(\frac{r_{4,s,m,p,q}^q}{1 - r_{4,s,m,p,q}^q}\right)^p
	\\&= 1.
\end{align*}
This shows the the radius $r_{4,s,m,p,q}$ is best possible.\vspace{1.2mm}

Similarly, for the operator $\mathscr{D}$, we can easily shown that  the inequality
\begin{align*}
	|\mathscr{D}f(\omega_m(z))|^s + |f(\omega_q(z))|^p&\le \frac{1}{2}+G_4(r).
\end{align*}
holds. The remaining part of the proof is analogous to the argument used for the inequality involving the operator $D$, and thus, we omit the details.
\end{proof}

\section{{\bf The Bohr phenomenon of the second-order harmonic differential operator.}}\label{Sec-4}
For \( \lambda \in [0,1] \), the convex combinations of a harmonic mapping and its differential operators $D$ and $\mathscr{D}$ are defined by  
\begin{align*}
	F_{\lambda}^{D} = (1 - \lambda) f + \lambda Df 
	\quad \text{and} \quad 
	F_{\lambda}^{\mathscr{D}} = (1 - \lambda) f + \lambda \mathscr{D}f,
\end{align*}
Furthermore, the series forms of the second-order harmonic differential operators, derived from equations (\ref{m1}) and (\ref{m2}), are given by
\begin{align*}
D^{2}f 
= D(Df) 
= z (Df)_{z} - \overline{z} (Df)_{\overline{z}} 
= z + \sum_{n=2}^{\infty} n^{2} a_{n} z^{n} 
    - \overline{\sum_{n=2}^{\infty} n^{2} b_{n} z^{n}}
\end{align*}
and 
\begin{align*}
	\mathscr{D}^{2}f 
	= \mathscr{D}(\mathscr{D}f) 
	= z (\mathscr{D}f)_{z} + \overline{z} (\mathscr{D}f)_{\overline{z}} 
	= z + \sum_{n=2}^{\infty} n^{2} a_{n} z^{n} 
	+ \overline{\sum_{n=2}^{\infty} n^{2} b_{n} z^{n}}.
\end{align*}
In the following result, we investigate an improved Bohr-type inequality for the linear combination of a univalent harmonic mapping and its associated differential operators. Considering $ f = h + \overline{g} \in \mathcal{S}_H^0 $, we establish sharp bounds involving the functions $ F_{\lambda}^{D} $ and $ F_{\lambda}^{\mathscr{D}} $ in terms of the Schwarz function $ \omega_m \in \mathcal{B}_m $.
\begin{theo}\label{Th-4.1}
Suppose $f = h + \overline{g} \in \mathcal{S}^0_H$ be a stable univalent harmonic mapping on the unit disk $\mathbb{D}$ and $\omega_m\in \mathcal{B}_m$ for all $m\geq 1$. Let $\lambda \in [0,1]$. Then 
\begin{align*}
|F_{\lambda}^{D}(\omega_m(z))|=|(1-\lambda)f(\omega_m(z))+\lambda D f(\omega_m(z))| \leq d\big(f(0), \partial f(\mathbb{D})\big) 
\end{align*}
and 
\begin{align*}
	|F_{\lambda}^{\mathscr{D}}(\omega_m(z))|=|(1-\lambda)f(\omega_m(z))+\lambda D f(\omega_m(z))| \leq d\big(f(0), \partial f(\mathbb{D})\big)
\end{align*}
both hold for $|z|=r \leq r_{5,m}$, where $r_{5,m}$ is the unique root in $(0,1)$ of the equation
\bea\label{rmp5}
\frac{(1-\lambda)r^m}{1-r^m}+\frac{\lambda r^m(1+r^m)}{(1-r^m)^3}-\frac{1}{4}=0\eea
The radius $r_{5,m}$ is the best possible.
\end{theo}
\begin{rem}
For the particular case $\omega_m(z) = z$ with $m= 1$, Theorem \ref{Th-4.1} reduces to \cite[Theorem 4.1]{Xu-An-Liu-JA-2025}. Hence, our theorem extends and improves the earlier result.
\end{rem}
\begin{proof}[{\bf Proof of Theorem \ref{Th-4.1}}]
Assume that $|z| = r$. By employing Lemma \ref{L1} together with Lemma \ref{L3}(i), we can estimate $|F_{\lambda}(\omega_m(z))|$ as follows:
\begin{align*}
|F_{\lambda}(\omega_m(z))|
&= \big|(1-\lambda)f(\omega_m(z)) + \lambda D f(\omega_m(z))\big| \\
&\le (1-\lambda)|f(\omega_m(z))| + \lambda |D f(\omega_m(z))| \\
&\le (1-\lambda)\frac{|\omega_m(z)|}{1-|\omega_m(z)|}
   + \lambda \bigg|\omega_m(z)
      + \sum_{n=2}^{\infty} n a_n (\omega_m(z))^n
      - \overline{\sum_{n=2}^{\infty} n b_n (\omega_m(z))^n}\bigg| \\
&\le (1-\lambda)\frac{|\omega_m(z)|}{1-|\omega_m(z)|}
   + \lambda \left(|\omega_m(z)| + \sum_{n=2}^{\infty} n(|a_n| + |b_n|)\,|\omega_m(z)|^n \right) \\
&\le (1-\lambda)\frac{|z|^m}{1-|z|^m}
   + \lambda \left(|z|^m + \sum_{n=2}^{\infty} n(|a_n| + |b_n|)\,|z|^{mn}\right) \\
&= (1-\lambda)\frac{r^m}{1-r^m}
   + \lambda \left(r^m + \sum_{n=2}^{\infty} n^2 r^{mn}\right) \\
&\le \frac{(1-\lambda)r^m}{1-r^m}
   + \frac{\lambda r^m(1+r^m)}{(1-r^m)^3} \\
&= \frac{1}{4} + G_5(r),
\end{align*}
where
\[
G_5(r)
   = \frac{(1-\lambda)r^m}{1-r^m}
     + \frac{\lambda r^m(1+r^m)}{(1-r^m)^3}
     - \frac{1}{4}.
\]

Clearly, $G_5(r)$ is continuous for all $r \in [0,1)$. In addition, we observe that
\[
G_5(0) = -\frac{1}{4}
\quad \text{and} \quad
\lim_{r \to 1^{-}} G_5(r) = +\infty.
\]
Differentiating $G_5(r)$ with respect to $r$ yields
\[
G_5'(r)
   = m r^{m-1}\,
     \frac{1 + (-2 + 6\lambda)r^m + r^{2m}}{(1-r^m)^4}
   \ge 0,
\]
which implies that $G_5(r)$ is non-decreasing on $[0,1)$.  
Consequently, $G_5(r)$ possesses a unique  zero in $(0,1)$.  
Denote this unique zero of equation~\eqref{rmp5} by $r_{5,m}$. Then for $r\geq r_{5,m}$,  we see that 
\begin{align}\label{T4.e} \frac{(1-\lambda)r^m}{1-r^m}+ \frac{\lambda r^m(1+r^m)}{(1-r^m)^3}\geq \frac{(1-\lambda)r_{5,m}^m}{1-r_{5,m}^m}+ \frac{\lambda r_{5,m}^m(1+r_{5,m}^m)}{(1-r_{5,m}^m)^3}
\end{align} 
Hence, the inequality
\[
|(1-\lambda)f(\omega_m(z))+\lambda D f(\omega_m(z))|
   \le d\big(f(0), \partial f(\mathbb{D})\big)
\]
holds for all $|z| \le r_{5,m}$.\vspace{1.2mm}

To show the radius $r_{5,m}$ is best possible. Henceforth, we consider  Koebe function $f=k$ given by \eqref{f1}, together with $\omega_m(z)=z^m$.  
For $r>r_{5,m}$, from (\ref{T4.e}), it follows that
\begin{align*}
|F_{\lambda}(r^m)|&=\frac{(1-\lambda)r^m}{1-r^m}+ \frac{\lambda r^m(1+r^m)}{(1-r^m)^3}\\
&> \frac{(1-\lambda)r_{5,m}^m}{1-r_{5,m}^m}
     + \frac{\lambda r_{5,m}^m(1+r_{5,m}^m)}{(1-r_{5,m}^m)^3}
   = \frac{1}{4}
   = d\big(f(0), \partial f(\mathbb{D})\big).
\end{align*}
Therefore, the radius $r_{5,m}$ cannot be improved, proving the result.\vspace{1.2mm}

Similarly, for the operator $\mathscr{D}$, we can easily show that the inequality$$|F_{\lambda}^{\mathscr{D}}(\omega_m(z))|=|(1-\lambda)f(\omega_m(z))+\lambda D f(\omega_m(z))| \le \frac{1}{2}+G_5(r)$$holds. Since the remaining part of the proof is analogous to the argument used for the inequality involving the operator $D$, we omit the details. This completes the proof.
\end{proof}
In the next result, we derive an improved Bohr-type inequality for the linear combination of a convex harmonic mapping and its corresponding differential operators. Considering $ f = h + \overline{g} \in \mathcal{S}_H^0 $, we obtain sharp estimates for the functions $ F_{\lambda}^{D} $ and $ F_{\lambda}^{\mathscr{D}} $ associated with the Schwarz function $ \omega_m \in \mathcal{B}_m $.
\begin{theo}\label{Th-4.2}
Suppose $f = h + \overline{g} \in \mathcal{S}^0_H$ be a stable convex harmonic mapping on the unit disk $\mathbb{D}$ and $\omega_m\in \mathcal{B}_m$ for all $m\geq 1$. Let $\lambda \in [0,1]$. Then 
\begin{align*}
|F_{\lambda}^{D}(\omega_m(z))|=|(1-\lambda)f(\omega_m(z))+\lambda D f(\omega_m(z))| \leq d\big(f(0), \partial f(\mathbb{D})\big) 
\end{align*}
and 
\begin{align}\label{Eq-4.3}
	|F_{\lambda}^{\mathscr{D}}(\omega_m(z))|=|(1-\lambda)f(\omega_m(z))+\lambda D f(\omega_m(z))| \leq d\big(f(0), \partial f(\mathbb{D})\big)
\end{align}
both hold for $|z|=r \leq r_{6,m}$, where $r_{6,m}$ is the unique  root in $(0,1)$ of the equation
\bea\label{rmp6}
\frac{(1-\lambda)r^m}{1-r^m}+\frac{\lambda r^m}{(1-r^m)^2}-\frac{1}{2}=0\eea
The radius $r_{6,m}$ is the best possible.
\end{theo}
\begin{rem}
Setting $\omega_m(z) = z$ with $m=p= 1$, Theorem \ref{Th-4.2} reduces to \cite[Corollary 4.2]{Xu-An-Liu-JA-2025}. Hence, our theorem extends and improves the earlier result.
\end{rem}
\begin{proof}
Let $|z| = r$. By applying Lemma \ref{L1} and Lemma \ref{L3}(i), we obtain the following estimate for $|F_{\lambda}(\omega_m(z))|$:
\begin{align*}
|F_{\lambda}(\omega_m(z))|
&= \big|(1-\lambda)f(\omega_m(z))+\lambda D f(\omega_m(z))\big| \\
&\le (1-\lambda)|f(\omega_m(z))| + \lambda |D f(\omega_m(z))| \\
&\le (1-\lambda)\frac{|\omega_m(z)|}{1-|\omega_m(z)|}
   + \lambda \bigg|\omega_m(z)
      + \sum_{n=2}^{\infty} n a_n (\omega_m(z))^n
      - \overline{\sum_{n=2}^{\infty} n b_n (\omega_m(z))^n}\bigg| \\
&\le (1-\lambda)\frac{|\omega_m(z)|}{1-|\omega_m(z)|}
   + \lambda \left(|\omega_m(z)| + \sum_{n=2}^{\infty} n(|a_n|+|b_n|)\,|\omega_m(z)|^n \right) \\
&\le (1-\lambda)\frac{|z|^m}{1-|z|^m}
   + \lambda \left(|z|^m + \sum_{n=2}^{\infty} n(|a_n|+|b_n|)\,|z|^{mn}\right) \\
&= (1-\lambda)\frac{r^m}{1-r^m}
   + \lambda \left(r^m + \sum_{n=2}^{\infty} n r^{mn}\right) \\
&\le \frac{(1-\lambda)r^m}{1-r^m}
   + \frac{\lambda r^m}{(1-r^m)^2} \\
&= \frac{1}{4} + G_6(r),
\end{align*}
where
\begin{align*}
G_6(r) = \frac{(1-\lambda)r^m}{1-r^m} + \frac{\lambda r^m}{(1-r^m)^2}- \frac{1}{2}.
\end{align*}
It is clear that $G_5(r)$ is continuous on $[0,1)$. Furthermore,
\begin{align*}
	G_6(0) = -\frac{1}{2}<0\;\text{and}\;\lim_{r \to 1^{-}} G_6(r) = +\infty.
\end{align*}
Differentiating $G_6(r)$ with respect to $r$, we find $G_6'(r)\geq 0$ for all $r\in [0, 1)$. Hence, $G_6(r)$ is a monotonically increasing function on $[0,1)$.  Consequently, $G_6(r)$ admits a unique  root in $(0,1)$.  
Let $r_{6,m}$ denote this unique root of equation \eqref{rmp6} in $(0,1)$. Then for $r\geq r_{6,m}$, we see that
\begin{align}\label{T4.e2} \frac{(1-\lambda)r^m}{1-r^m}+ \frac{\lambda r^m}{(1-r^m)^2}\geq \frac{(1-\lambda)r_{6,m}^m}{1-r_{6,m}^m}+ \frac{\lambda r_{6,m}^m}{(1-r_{6,m}^m)^2}.\end{align}
 
Therefore, the desired inequality
\[
|(1-\lambda)f(\omega_m(z))+\lambda D f(\omega_m(z))|
   \le d\big(f(0), \partial f(\mathbb{D})\big)
\]
holds for $|z| \le r_{6,m}$.\vspace{1.2mm}

To establish radius $r_{6,m}$ is the best possible, we consider the function $f_2$ defined in \eqref{f2} together with $\omega_m(z)=z^m$.  
For $r>r_{6,m}$, from (\ref{T4.e2}) we compute
\begin{align*}
|F_{\lambda}(r^m)|&=\frac{(1-\lambda)r^m}{1-r^m}+ \frac{\lambda r^m}{(1-r^m)^2}\\
 &> \frac{(1-\lambda)r_{6,m}^m}{1-r_{6,m}^m}
     + \frac{\lambda r_{6,m}^m}{(1-r_{6,m}^m)^2}
   = \frac{1}{2}
   = d\big(f_2(0), \partial f_2(\mathbb{D})\big).
\end{align*}
Thus, the radius $r_{6,m}$ is the best possible.\vspace{1.2mm}

Similarly, for the operator $\mathscr{D}$, we can easily show that the inequality \eqref{Eq-4.3} holds. The remaining part of the proof is analogous to the argument used for the inequality involving the operator $D$, and thus, we omit the details. This completes the proof.
\end{proof}
We now present an improved Bohr-type inequality for the second-order differential operators of stable univalent harmonic mappings. For $ f = h + \overline{g} \in \mathcal{S}_H^0 $, we obtain sharp inequalities involving $ D^2f $ and $ \mathscr{D}^2f $ corresponding to the Schwarz functions $ \omega_m \in \mathcal{B}_m $ and $ \omega_p \in \mathcal{B}_p $.
\begin{theo}\label{Th-4.3}
Suppose $f = h + \overline{g} \in \mathcal{S}^0_H$ be a stable univalent harmonic mapping on the unit disk $\mathbb{D}$ and $\omega_m\in \mathcal{B}_m$ for all $m\geq 1$. Then 
\begin{align*}
|D^2 f(\omega_m(z))|+\sum_{n=2}^{\infty} (|a_n|+|b_n|)|\omega_p(z)|^n \leq d\big(f(0), \partial f(\mathbb{D})\big)
\end{align*}
and 
\begin{align}\label{Eq-4.6}
	|\mathscr{D}^2 f(\omega_m(z))|+\sum_{n=2}^{\infty} (|a_n|+|b_n|)|\omega_p(z)|^n \leq d\big(f(0), \partial f(\mathbb{D})\big)
\end{align}
both hold for $|z|=r \leq r_{7,m,p}$, where $r_{7,m,p}$ is the unique root in $(0,1)$ of the equation
\bea\label{rmp7}
\frac{r^m(1+4r^m+r^{2m})}{(1-r^m)^4}+\frac{r^{pN}(N+r^p-Nr^p)}{(1-r^p)^2}-\frac{1}{4}=0\eea
The radius $r_{6,m}$ is the best possible.
\end{theo}
\begin{rem}
Setting $\omega_m(z) = z$, $\omega_p(z) = z$ with $m=p= 1$, Theorem \ref{Th-4.3} reduces to \cite[Theorem 4.5]{Xu-An-Liu-JA-2025}. Hence, our theorem extends and improves the earlier result.
\end{rem}
\begin{proof}
Let $|z| = r$. By applying Lemma~\ref{L1} and Lemma~\ref{L3}(i), we obtain the following estimate for $|F_{\lambda}(\omega_m(z))|$:
\begin{align*}
&|D^2 f(\omega_m(z))|+\sum_{n=N}^{\infty} (|a_n|+|b_n|)|\omega_p(z)|^n\\
&\le \bigg|\omega_m(z)+\sum_{n=2}^{\infty} n^2 a_n (\omega_m(z))^n
      - \overline{\sum_{n=2}^{\infty} n^2 b_n (\omega_m(z))^n}\bigg| + \sum_{n=N}^{\infty} (|a_n|+|b_n|)|z|^{np} \\
&\quad+|\omega_m(z)|+\sum_{n=2}^{\infty} n^2(|a_n|+|b_n|)|\omega_m(z)|^n +\sum_{n=N}^{\infty} (|a_n|+|b_n|)|z|^{np} \\
&\le |z|^m+ \sum_{n=2}^{\infty} n^2(|a_n|+|b_n|)|z|^{mn}+\sum_{n=N}^{\infty} (|a_n|+|b_n|)|z|^{np}\\
&\le r^m+ \sum_{n=2}^{\infty} n^3r^{mn}+\sum_{n=N}^{\infty}nr^{np}\\
&\le \frac{r^m(1+4r^m+r^{2m})}{(1-r^m)^4}+\frac{r^{pN}(N+r^p-Nr^p)}{(1-r^p)^2} \\
&= \frac{1}{4} + G_7(r),
\end{align*}
where
\[
G_7(r) = H_1(r^m)+H_2(r^p) -\frac{1}{4}
\]
with 
\begin{align*}
	H_1(r^m)=\frac{r^m(1+4r^m+r^{2m})}{(1-r^m)^4}\; \mbox{and}\;H_2(r^m)=\frac{r^{pN}(N+r^p-Nr^p)}{(1-r^p)^2}.
\end{align*} Setting $t = r^m$ and $s = r^p$, we obtain
\begin{align*}
	H_1(t)= \frac{t(1+4t+t^{2})}{(1-t)^4}\; \mbox{and}\;H_2(s)=\frac{s^{N}(N+s-Ns)}{(1-s)^2}
\end{align*} 
It is clear that $G_7(r)$ is continuous on $[0,1)$. Furthermore,
\begin{align*}
	G_7(0) = -\frac{1}{4}<0
	\quad \text{and} \quad
	\lim_{r \to 1^{-}} G_7(r) = +\infty.
\end{align*}
Differentiating $G_7(r)$ with respect to $r$, we obtain
\[
G_7'(r) = mH_1'(r^m)r^{m-1}+ pH_2'(r^p)r^{p-1}.
\]
It is easy to see that
\[H_1'(t) = \frac{1 + 11t + 11t^2 + t^3}{(1 - t)^5}\geq 0\]
 and 
\beas
H_2'(s)=\frac{s^{N-1}\big[N^2 + s(1 - N^2)\big](1 - s)
      + 2s^N(N + s - Ns)}{(1 - s)^3}\geq 0.
\eeas
Thus, it is clear that $G_7'(r) \geq 0$, which means $G_7(r)$ is a monotonically increasing function on $[0,1)$. Consequently, $G_7(r)$ admits a unique root in $(0,1)$. Denoting this unique root of \eqref{rmp7} as $r_{7,m,p}$, we find that for $r \geq r_{7,m,p}$, 
\begin{align}\label{T4.e3}&\frac{r^m(1+4r^m+r^{2m})}{(1-r^m)^4}+\frac{r^{pN}(N+r^p-Nr^p)}{(1-r^p)^2}\nonumber\\
  &\geq \frac{r_{7,m,p}^m(1+4r_{7,m,p}^m+r_{7,m,p}^{2m})}{(1-r_{7,m,p}^m)^4}+\frac{r_{7,m,p}^{pN}(N+r_{7,m,p}^p-Nr_{7,m,p}^p)}{(1-r_{7,m,p}^p)^2}.
\end{align} 
Thus, the desired inequality
\begin{align*}
	|D^2 f(\omega_m(z))|+\sum_{n=N}^{\infty} (|a_n|+|b_n|)|\omega_p(z)|^n
	\le d\big(f(0), \partial f(\mathbb{D})\big)
\end{align*}
is obtained for $|z| \le r_{7,m,p}$.\vspace{1.2mm}

The next part of the proof is to show that the radius $r_{7,m,p}$ is the best possible. We consider the function $f=k$ defined in Equation~\eqref{f1} together with $\omega_m(z)=z^m$. For $r>r_{7,m,p}$, in view of (\ref{T4.e3}), a simple computation shows that
\begin{align*}
&|D^2 f_1(r^m)|+\sum_{n=N}^{\infty} (|a_n|+|b_n|)|r^p|^n\\
&=\frac{r^m(1+4r^m+r^{2m})}{(1-r^m)^4}+\frac{r^{pN}(N+r^p-Nr^p)}{(1-r^p)^2}\\
 &  > \frac{r_{7,m,p}^m(1+4r_{7,m,p}^m+r_{7,m,p}^{2m})}{(1-r_{7,m,p}^m)^4}+\frac{r_{7,m,p}^{pN}(N+r_{7,m,p}^p-Nr_{7,m,p}^p)}{(1-r_{7,m,p}^p)^2}\\
 &= \frac{1}{4}\\&= d\big(f(0), \partial f(\mathbb{D})\big).
\end{align*}
Thus, the radius $r_{7,m,p}$ is the best possible.\vspace{1.2mm}

Similarly, the inequality \eqref{Eq-4.6} is easily shown to hold for the operator $\mathscr{D}$. Since the rest of the proof is analogous to the argument for the operator $D$, we omit the details. This completes the proof.
\end{proof}

Next, we  derive an improved Bohr-type inequality for the second-order differential operators of stable convex harmonic mappings. For $ f = h + \overline{g} \in \mathcal{S}_H^0 $, we obtain sharp inequalities involving $ D^2f $ and $ \mathscr{D}^2f $ corresponding to the Schwarz functions $ \omega_m \in \mathcal{B}_m $ and $ \omega_p \in \mathcal{B}_p $.
\begin{theo}\label{Th-4.4}
Suppose $f = h + \overline{g} \in \mathcal{S}^0_H$ be a stable convex harmonic mapping on the unit disk $\mathbb{D}$ and $\omega_m\in \mathcal{B}_m$, $\omega_p\in \mathcal{B}_p$. Then 
\begin{align*}
|D^2 f(\omega_m(z))|+\sum_{n=2}^{\infty} (|a_n|+|b_n|)|\omega_p(z)|^n \leq d\big(f(0), \partial f(\mathbb{D})\big)
\end{align*}
and 
\begin{align}\label{Eq-4.9}
	|\mathscr{D}^2 f(\omega_m(z))|+\sum_{n=2}^{\infty} (|a_n|+|b_n|)|\omega_p(z)|^n \leq d\big(f(0), \partial f(\mathbb{D})\big)
\end{align}
both hold for $|z|=r \leq r_{8,m,p}$, where $r_{8,m,p}$ is the unique  root in $(0,1)$ of the equation
\bea\label{rmp8}
 \frac{r^m(1 + r^m)}{(1 - r^m)^3} + \frac{r^{pN}}{1 - r^p} - \frac{1}{2}=0\eea
The radius $r_{8,m,p}$ is the best possible.
\end{theo}
\begin{rem}
Setting $\omega_m(z) = z$, $\omega_p(z) = z$ with $m=p= 1$, Theorem \ref{Th-4.4} reduces to \cite[Corollary 4.6]{Xu-An-Liu-JA-2025}. Hence, our theorem improves the earlier result.
\end{rem}
\begin{proof}
Let $|z| = r$. By applying Lemma~\ref{L1} and Lemma~\ref{L3}(i), we obtain the following estimate:
\begin{align*}
&|D^2 f(\omega_m(z))| + \sum_{n=N}^{\infty} (|a_n| + |b_n|)\,|\omega_p(z)|^n \\
&\le \bigg|\omega_m(z) + \sum_{n=2}^{\infty} n^2 a_n \big(\omega_m(z)\big)^n
      - \overline{\sum_{n=2}^{\infty} n^2 b_n \big(\omega_m(z)\big)^n}\bigg|
      + \sum_{n=N}^{\infty} (|a_n| + |b_n|)\,|z|^{np} \\
&\le |\omega_m(z)| + \sum_{n=2}^{\infty} n^2(|a_n| + |b_n|)\,|\omega_m(z)|^n
      + \sum_{n=N}^{\infty} (|a_n| + |b_n|)\,|z|^{np} \\
&\le |z|^m + \sum_{n=2}^{\infty} n^2(|a_n| + |b_n|)\,|z|^{mn}
      + \sum_{n=N}^{\infty} (|a_n| + |b_n|)\,|z|^{np} \\
&\le r^m + \sum_{n=2}^{\infty} n^2 r^{mn} + \sum_{n=N}^{\infty} r^{np} \\
&\le \frac{r^m(1 + r^m)}{(1 - r^m)^3} + \frac{r^{pN}}{1 - r^p} \\
&= \frac{1}{2} + G_8(r),
\end{align*}
where
\[
G_8(r) = \frac{r^m(1 + r^m)}{(1 - r^m)^3} + \frac{r^{pN}}{1 - r^p} - \frac{1}{2}.
\]

It is clear that $G_8(r)$ is continuous on $[0,1)$. Furthermore,
\[
G_8(0) = -\frac{1}{2}<0
\quad \text{and} \quad
\lim_{r \to 1^{-}} G_8(r) = +\infty.
\]
Differentiating $G_8(r)$ with respect to $r$, we obtain
\[
G_8'(r)
= \frac{m r^{m-1}\big(1 + 4r^m + r^{2m}\big)}{(1 - r^m)^4}
  + \frac{p r^{pN-1}\big(N + (1 - N)r^p\big)}{(1 - r^p)^2} \ge 0.
\]
Hence, $G_8(r)$ is a monotonically increasing function on $[0,1)$.  
Consequently, $G_8(r)$ admits a unique  root in $(0,1)$.  
Let $r_{8,m,p}$ denote this unique root of equation~\eqref{rmp7} in $(0,1)$. For $r\geq r_{8,m,p}$, we see that
\begin{align}\label{T4.e4} \frac{r^m(1 + r^m)}{(1 - r^m)^3} + \frac{r^{pN}}{1 - r^p}\geq \frac{r_{8,m,p}^m(1 + r_{8,m,p}^m)}{(1 - r_{8,m,p}^m)^3}
     + \frac{r_{8,m,p}^{pN}}{1 - r_{8,m,p}^p}.
\end{align} 
 Thus, we see that
\begin{align*}
	|D^2 f(\omega_m(z))| + \sum_{n=N}^{\infty} (|a_n| + |b_n|)\,|\omega_p(z)|^n
	\le d\big(f(0), \partial f(\mathbb{D})\big)
\end{align*}
holds for $|z| \le r_{8,m,p}$.\vspace{1.2mm}

To establish that the radius $r_{8,m,p}$ is the best possible, we consider the function $f_2$ defined in Equation~\eqref{f2} together with $\omega_m(z) = z^m$. For $r > r_{8,m,p}$, from (\ref{T4.e4}), a simple computation leads to
\begin{align*}
|D^2 f_2(r^m)| &+ \sum_{n=N}^{\infty} (|a_n| + |b_n|)\,|r^p|^n\\
&=\frac{r^m(1 + r^m)}{(1 - r^m)^3} + \frac{r^{pN}}{1 - r^p}\\
  &> \frac{r_{8,m,p}^m(1 + r_{8,m,p}^m)}{(1 - r_{8,m,p}^m)^3}
     + \frac{r_{8,m,p}^{pN}}{1 - r_{8,m,p}^p} \\
  &= \frac{1}{2}\\&
   = d\big(f_2(0), \partial f_2(\mathbb{D})\big).
\end{align*}
Thus, the radius $r_{8,m,p}$ is the best possible.\vspace{1.2mm}

Similarly, for the operator $\mathscr{D}$, we can easily show that the inequality \eqref{Eq-4.9} holds. The remaining part of the proof is analogous to the argument used for the inequality involving the operator $D$; thus, we omit the details here. This completes the proof.
\end{proof}
\section{\bf Improved Bohr inequality in terms of $S_r/\pi$ for a class of stable harmonic mappings}\label{Sec-5}
  
 The area $S_r$ of the disk $|z| < r$ under the harmonic map $f_d$ is given by  
\begin{align*}
	\frac{S_r}{\pi} = \frac{1}{\pi} \iint_{D_r} J_{f_d} \, dx \, dy 
	= \frac{1}{\pi} \iint_{D_r} \left( |h_d'(z)|^2 - |g_d'(z)|^2 \right) dx \, dy.
\end{align*}
Because $f_d$ is sense-preserving, we have 
\begin{align*}
	J_{f_d}(z) = |h_d'(z)|^2 - |g_d'(z)|^2 > 0.
\end{align*}
In the study of improved Bohr inequalities for various function classes, the quantity $S_r/\pi$ plays a significant role. However, the Bohr inequality involving Schwarz functions $\omega_m \in \mathcal{B}_m$ together with a polynomial expression of the quantity $S_r/\pi$ has not yet been examined in the literature. Given the importance of the Bohr inequality and its sharpness for the class $\mathcal{S}_H^0$, it is therefore natural to pose the following question:
\begin{ques}\label{Qn-5.1}
	Can we establish an improved Bohr inequality involving Schwarz functions $\omega_m \in \mathcal{B}_m$ and a polynomial expression of the quantity $S_r/\pi$, for the class of functions $f=h+\overline{g} \in \mathcal{S}_H^0$? Furthermore, can we prove that the resulting radius is best possible (sharp)?
\end{ques}
In this section, we prove the following result, which affirmatively answers Question \ref{Qn-5.1}.
\begin{theo}\label{Th-5.1} Let $f=h+\overline  g\in \mathcal{S}_H^0$ and let $P_k(t)=\sum_{i=1}^\infty a_i t^i$, where $a_i\geq 0$. Then 
	\begin{align*}
		M(|\omega_m(z)|)+P_k\bigg(\frac{S_r}{\pi}\bigg)\leq d(f(0),\partial f(\mathbb{D}))\; \mbox{for}\; |z|\leq r_{m,k},
	\end{align*}
	 where $r_{m,k}$ is the unique root in $(0,1)$ of the equation
	 \begin{align*}
	 	\frac{r^m}{(1-r^m)^2}+P_k\left(\frac{r^2+26r^4+66r^6+26r^8+r^{10}}{(1-r^2)^6}\right)-\frac{1}{4}=0.
	 \end{align*}
	 The radius $r_{m,k}$ is best possible.
\end{theo} 
\begin{rem} The following observations are clear.
	\begin{enumerate}
		\item[(i)] In particular, when we set $m=1$ (so $\omega_1 \in \mathcal{B}_1$) and choose the polynomial $P_k(t)=t$, Theorem \ref{Th-5.1} simplifies to \cite[Theorem 3.7]{Xu-An-Liu-JA-2025}, demonstrating a significant improvement of that result.
		\item[(ii)] In fact, for  $m=1$ (so $\omega_1 \in \mathcal{B}_1$), $a_1=16/9$ and $a_2=\lambda=18.6095...$ and $F(t)=a_3t^3+\cdots+a_kt^k$, $(k\geq 3)$, we see that $F : [0, \infty)\to [0, \infty)$ and in this case, our inequality becomes
		\begin{align}\label{Eq-5.1}
			M(|z|)+\frac{16}{9}\left(\frac{S_r}{\pi}\right)+\lambda\left(\frac{S_r}{\pi}\right)^2+F\bigg(\frac{S_r}{\pi}\bigg)\leq d(f(0),\partial f(\mathbb{D}))
		\end{align}
		for $|z|\leq r_{1,k},$ where $\leq r_{1,k}$ is the unique root in $(0, 1)$ of the equation 
		\begin{align*}
			&\frac{16}{9}\left(\frac{r^2+26r^4+66r^6+26r^8+r^{10}}{(1-r^2)^6}\right)+\lambda\left(\frac{r^2+26r^4+66r^6+26r^8+r^{10}}{(1-r^2)^6}\right)^2\\&\quad+\frac{r}{(1-r)^2}+F\left(\frac{r^2+26r^4+66r^6+26r^8+r^{10}}{(1-r^2)^6}\right)-\frac{1}{4}=0.
		\end{align*}
		The radius $r_{1,k}$ is best possible. Our Theorem \ref{Th-5.1} for harmonic mappings overcomes the limitation pointed out in \cite[Remark 2]{Ismagilov-2020-JMAA}, which stated that inequality \eqref{Eq-5.1} failed to hold for the class $\mathcal{B}$ of analytic functions.
	\end{enumerate}
\end{rem}
\begin{proof}
	Switching to polar coordinates and using Parseval's identity gives
	\begin{align*}
		\frac{S_r}{\pi}
		&= \sum_{n=1}^{\infty} n^3 \left( |a_n|^2 - |b_n|^2 \right) r^{2n}\\
		&= \sum_{n=1}^{\infty} n^3 (|a_n| - |b_n|)(|a_n| + |b_n|) r^{2n}\\
		&= \sum_{n=1}^{\infty} n^3 (\, ||a_n| - |b_n||\,)(|a_n| + |b_n|) r^{2n}.
	\end{align*}

	By applying Lemma \ref{L1}, we obtain the following estimate
	\begin{align} M(|\omega_m(z)|)+P_k\bigg(\frac{S_r}{\pi}\bigg)&\leq |\omega_m(z)|+\sum_{n=2}^{\infty} |a_n||\omega_m(z)|^n+\sum_{n=2}^{\infty} |b_n||\omega_m(z)|^n+P_k\left(\sum_{n=1}^\infty n^5r^{2n}\right)\nonumber\\
		&\leq |z|^m+\sum_{n=2}^{\infty} (|a_n|+|b_n|)|z|^{mn}+P_k\left(\sum_{n=1}^\infty n^5r^{2n}\right)\nonumber\\
		& \leq r^m+\sum_{n=2}^{\infty}nr^{mn}+P_k\left(\frac{r^2+26r^4+66r^6+26r^8+r^{10}}{(1-r^2)^6}\right)\nonumber\\
		&=\frac{r^m}{(1-r^m)^2}+P_k\left(\frac{r^2+26r^4+66r^6+26r^8+r^{10}}{(1-r^2)^6}\right)\nonumber\\
		&=\frac{1}{4}+G(r),
	\end{align}
	where 
	\[G(r)=\frac{r^m}{(1-r^m)^2}+P_k\left(\frac{r^2+26r^4+66r^6+26r^8+r^{10}}{(1-r^2)^6}\right)-\frac{1}{4}.\]
	It follows directly from the definition that $G(r)$ is differentiable on $[0,1)$. Differentiating $G(r)$ with respect to $r$, we obtain 
	\[ G'(r)=\frac{mr^{m-1}}{(1-r^m)^2}+\frac{2r^m}{(1-r^m)^3}mr^{m-1}+P(r),\]
	where 
	\begin{align}\label{jk1} P(r)=&\bigg(\sum_{i=1}^{k-1}ia_i\left(\frac{r^2+26r^4+66r^6+26r^8+r^{10}}{(1-r^2)^6}\right)^{i-1}\bigg)H_0(r),
	\end{align}
	where \[H_0(r):=\frac{2r+104r^3+369r^5+208r^7+10r^9}{(1-r^2)^6}+\frac{12(r^2+26r^4+66r^6+26r^8+r^{10})r}{(1-r^2)^6)}.\]
	Because all terms in $G'(r)$ are positive, the function $G(r)$ is strictly increasing on the interval $[0, 1)$. Consequently, for all $r \geq r_{m,k}$, the inequality $G(r) \geq G(r_{m,k})$ holds, implying the result
	\begin{align}\label{Eq-5.4}
		\frac{r^m}{(1-r^m)^2}&+P_k\left(\frac{r^2+26r^4+66r^6+26r^8+r^{10}}{(1-r^2)^6}\right)\\
		&\geq \frac{r_{m,k}^m}{(1-r_{m,k}^m)^2}+P_k\left(\frac{r_{m,k}^2+26r_{m,k}^4+66r_{m,k}^6+26r_{m,k}^8+r^{10}}{(1-r_{m,k}^2)^6}\right).\nonumber
	\end{align}
	Moreover, it is easy to see that 
	\begin{align*}
		G(0)=-\frac{1}{4}<0\;\;\;\text{and}\;\;\;\lim_{r\rightarrow 1^{-1}}G(r)=+\infty.
	\end{align*}
	Consequently, $G(r)$ has a unique root in the interval $(0, 1)$. Denoting this unique root by $r_{m,k}$, we obtain the desired inequality
	\begin{align*}
		M(|\omega_m(z)|)+P_k\bigg(\frac{S_r}{\pi}\bigg)\leq d(f(0),\partial f(\mathbb{D}))\; \mbox{for}\; |z|=r\leq r_{m,k}.
	\end{align*}
	
	To show the radius $r_{m,k}$ is best possible, we consider the Koebe function 
	\begin{align*}
		k(z)=\frac{z}{(1-z)^2}=z+\sum_{n-2}^{\infty}nz^n,
	\end{align*}
	which gives $a_n=n$ and $b_n=0$.\vspace{1.2mm}
	
	A simple computation shows for the function $k(z)$ that
	\begin{align*}\frac{S_r}{\pi}&=\sum_{n=1}^{\infty} n^3(||a_n|-|b_n||(|a_n|+|b_n|)r^{2n}\\
		&=\sum_{n=1}^{\infty} n^5r^{2n}=\frac{r^2+26r^4+66r^6+26r^8+r^{10}}{(1-r^2)^6}.\end{align*}
	For $\omega_m(z)=z^m (m\geq 1)$ and $z=r>r_{m,k}$, with $f=k$, it follows from (\ref{jk1}), \eqref{Eq-5.4} and \eqref{Eq-3.4} that
	\begin{align*}
		M(|\omega_m(z)|)+P_k\bigg(\frac{S_r}{\pi}\bigg)=& M(|z^m|)+P_k\bigg(\frac{S_r}{\pi}\bigg)\\
		&= r^m+\sum_{n=2}^{\infty}nr^{mn}+P_k\left(\frac{r^2+26r^4+66r^6+26r^8+r^{10}}{(1-r^2)^6}\right)\\
		&=\frac{r^m}{(1-r^m)^2}+P_k\left(\frac{r^2+26r^4+66r^6+26r^8+r^{10}}{(1-r^2)^6}\right)\\
		&> \frac{r_{m,k}^m}{(1-r_{m,k}^m)^2}+P_k\left(\frac{r_{m,k}^2+26r_{m,k}^4+66r_{m,k}^6+26r_{m,k}^8+r^{10}}{(1-r_{m,k}^2)^6}\right)\\
		&=\frac{1}{4}\\
		&=d(k(0),\partial k(\mathbb{D}).
	\end{align*}  
	which shows that the radius $r_{m,k}$ is the best possible. 
\end{proof}
\noindent{\bf Acknowledgment:} The authors would like to thank the referees for their suggestions and comments to improve exposition of the paper. \vspace{1.2mm}

\noindent {\bf Funding:} Not Applicable.\vspace{1.2mm}

\noindent\textbf{Conflict of interest:} The authors declare that there is no conflict  of interest regarding the publication of this paper.\vspace{1.2mm}

\noindent\textbf{Data availability statement:}  Data sharing not applicable to this article as no datasets were generated or analysed during the current study.\vspace{1.2mm}

\noindent {\bf Authors' contributions:} All the authors have equal contributions in preparation of the manuscript.

\end{document}